\documentclass[11pt,a4paper]{article}
\usepackage{amsthm,amsfonts,amsmath,amscd}

\title{Semi-infinite $q$-wedge construction of
  \\
  the level~2 Fock Space of $U_q(\affsl{2})$ \thanks{This proceedings
    is based on the paper ``Perfect crystals and $q$-Fock spaces''
    \cite{KMPY} by Masaki Kashiwara, Tetsuji Miwa, Chung Ming Yung and
    the speaker.}}

\author{Jens-Ulrik Holger Petersen \thanks{extended version
    of an invited talk at the Nankai-CRM meeting ``Extended and
    quantum algebras and their application to physics'', held at the
    Nankai Institute for Mathematics, Nankai University, Tianjin,
    China, 19--24 August 1996}}

\date{January 1997\\
  q-alg/9701040}

\newcommand{\N}{{\mathbb{N}}}
\newcommand{\Q}{{\mathbb{Q}}}
\newcommand{\Z}{{\mathbb{Z}}}
\newcommand{\Zplus}{\Z_{>0}}
\newcommand{\brackets}[1]{\left[ #1 \right]} 
\newcommand{\range}[1]{\brackets{#1}} 

\renewcommand{\b}[1]{b_{#1}}
\newcommand{\bo}[1]{\b{#1}^\circ}
\renewcommand{\v}[1]{v_{#1}}
\newcommand{\vo}[1]{\v{#1}^\circ}
\newcommand{\tensor}{\otimes}

\newcommand{\Directsum}{\bigoplus}
\renewcommand{\sl}[1]{{\mathfrak{sl}_{#1}}}
\newcommand{\affsl}[1]{\hat{\sl{}}_{#1}} 
\newcommand{\g}{\mathfrak{g}}
\newcommand{\pairing}[1]{\langle{#1}\rangle}
\newcommand{\Vaff}{V_{\aff}}
\newcommand{\Laff}{L_{\aff}}
\newcommand{\Baff}{B_{\aff}}

\newcommand{\Pcl}{P_{\cl}}
\newcommand{\F}[1]{\mathcal{F}_{#1}}

\newcommand{\set}[1]{\{ #1 \}}
\newcommand{\vac}[1]{|{#1}\rangle}

\newcommand{\U}{{\mathrm{U}}\hspace{-1.5pt}}

\newcommand{\Wedge}{{\textstyle\bigwedge}\!}
\newcommand{\te}{\tilde{e}}
\newcommand{\tf}{\tilde{f}}
\newcommand{\e}{\varepsilon}
\newcommand{\f}{\varphi}
\newcommand{\card}[1]{|#1|}
\renewcommand{\H}{\mathrm{H}}

\DeclareMathSymbol{\rightleftarrows}{\mathrel}{AMSa}{"1D}
\DeclareMathOperator*{\rightleftcrystal}{\rightleftarrows}
\DeclareMathOperator{\wt}{wt}
\DeclareMathOperator{\cl}{cl}
\DeclareMathOperator{\aff}{aff}
\DeclareMathOperator{\Span}{span}
\DeclareMathOperator{\ch}{ch}

\theoremstyle{plain}
\newtheorem{lemma}{Lemma}[section]
\newtheorem{proposition}[lemma]{Proposition}

\newtheorem{theorem}[lemma]{Theorem}

\theoremstyle{definition}

\theoremstyle{remark}

\numberwithin{equation}{subsection}

\begin{document}
\maketitle

%
\vspace{-2em}
\begin{abstract}
  In this proceedings a particular example from~\cite{KMPY} is
  presented: the construction of the level~2 Fock space of
  $\U_q(\affsl{2})$.  The generating ideal of the wedge relations is
  given and the wedge space defined.  Normal ordering of wedges is
  defined in terms of the energy function.  Normally ordered wedges
  form a base of the wedge space.

  The $q$-deformed Fock space is defined as the space of semi-infinite
  wedges with a finite number of vectors in the wedge product
  differing from a ground state sequence and endowed with a separated
  $q$-adic topology .  Normally ordered wedges form a base of the Fock
  space.  The action of $\U_q(\affsl{2})$ on the Fock space converges
  in the $q$-adic topology.  On the Fock space the action of bosons,
  which commute with the $\U_q(\affsl{2})$-action, also converges in
  the $q$-adic topology.  Hence follows the decomposition of the Fock
  space into irreducible $\U_q(\affsl{2})$-modules.
\end{abstract}

\section{Introduction}
The classical semi-infinite wedge construction~\cite{DJKM} originates
from the study of the representation theory of affine Lie algebras and
the soliton theory of integrable hierarchies during the 80's (for a
review see for example~\cite{KR}).  In the
last couple of years the subject has been revived by the consideration
of its $q$-deformation in the context of the representation theory of
quantum affine algebras.

In~\cite{S,KMS} the semi-infinite wedge space construction of the
level~1 Fock space of $\U_q(\affsl{n})$ and its decomposition were
given.  In~\cite{KMPY} under certain assumptions a general scheme for
the wedge construction of $q$-deformed Fock spaces using the theory of
perfect crystals was presented.
Let $\U_q(\g)$ be a quantum affine algebra.  Let $V$ be a finite
dimensional $\U_q(\g')$-module with a perfect crystal base~\cite{KMN}
of level~$k$.  Let $\Vaff$ denote its affinization.
In~\cite{KMPY} the wedge space $\Wedge^r\Vaff$ ($r\in\N$) was
constructed using the $\U_q(\g)$-action and the following
$\U_q(\g)$-linear map was given
\begin{displaymath}
  \Wedge^r\Vaff\tensor\Wedge^s\Vaff \longrightarrow \Wedge^{r+s}\Vaff.
\end{displaymath}
Normal ordering of wedges was defined and it was proven that normally
ordered wedges form a base of $\Wedge^r\Vaff$.
The level~$k$ Fock space $\F{m}$ ($m\in\Z$) was constructed and the
corresponding $\U_q(\g)$-linear map
\begin{displaymath}
  \Wedge^r\Vaff\tensor \F{m} \longrightarrow \F{m-r}
\end{displaymath}
was defined.  The decomposition of the Fock space was given.

In~\cite{KMPY} examples of the theory for level~1 $A^{(1)}_n$
(\cite{KMS}), $B^{(1)}_n$, $D^{(1)}_n$, $A^{(2)}_{2n}$,
$A^{(2)}_{2n-1}$, $D^{(2)}_{n+1}$ and level~$k$ $A^{(1)}_1$ were also
given.

In this talk I summarize some of the results of~\cite{KMPY},
describing the case of level~2 $\U_q(\affsl{2})$ in detail.  This
fairly simple example is sufficient to illustrate the differences
between the level~1 $\U_q(\affsl{n})$ case (\cite{KMS}) and other
cases (\cite{KMPY}): in particular the need to endow the Fock space
with a separated $q$-adic topology.  For fuller details and proofs,
please see~\cite{KMPY}.

\subsection{Thanks}
I would like to thank the organisers --- Professor Mo-Lin Ge,
Professor Yvan Saint-Aubin and Professor Luc Vinet --- for the
invitation and a most enjoyable meeting at Nankai University.  I thank
Professor Ge and all the local organisers for the excellent
hospitality.  I thank R.I.M.S. for travel support.  Finally I thank my
coauthors of~\cite{KMPY} for the collaboration, which made this paper
possible.

\section{Preliminaries}
Let $\g$ be an affine Lie algebra with associated weight lattice
$P:=\sum_{i\in I}\Z\Lambda_i\oplus\Z\delta$ and $\g'$ its derived Lie
subalgebra with associated weight lattice $\Pcl:=\sum_{i\in
  I}\Z\Lambda^{\cl}_i$.  Define $h_i\in P^*$ by
$\pairing{h_i,\delta}=0$ and $\pairing{h_i,\Lambda_j}=\delta_{i,j}$
($i,j\in I$).  Similarly define $h_i\in\Pcl^*$ by
$\pairing{h_i,\Lambda^{\cl}_j}=\delta_{i,j}$ ($i,j\in I$).  Let $W$ be
the Weyl group of $\g$.  Let $c$ denote the canonical central element
in $\g$ and $\g'$.  Define the projection $\cl:P \rightarrow\Pcl$ by
\begin{displaymath}
  \cl\bigl(\sum_{i\in I} \omega_i\Lambda_i +\omega_\delta\delta\bigr) :=
  \sum_{i\in I} \omega_i\Lambda_i^{\cl} \qquad
  (\omega_i,\omega_\delta\in\Z).
\end{displaymath}
Define $\Pcl^+:=\set{\lambda\in\Pcl\mid \pairing{h_i,\lambda}\geq 0}$.

Let $\U_q(\g)$ (respectively $\U_q(\g')$) be the quantum universal
enveloping algebra of $\g$~($\g'$) over $\Q(q)$ with generators
$e_i,f_i,q^h$ ($i\in I$ and $h\in P^*$ ($\Pcl^*$)).  Write $t_i$ for
$q_i^{h_i}$.

The coproduct of $\U_q(\g)$ $(\U_q(\g'))$ is taken to be
$\Delta=\bar{\Delta}_+$:
\begin{displaymath}
  \Delta:
  \begin{cases}
    q^h\mapsto q^h\tensor q^h &(h\in P^*\: (\Pcl^*))
    \\
    e_i\mapsto e_i\tensor 1 +t_i^{-1}\tensor e_i &(i\in I)
    \\
    f_i\mapsto f_i\tensor t_i +1\tensor f_i &(i\in I)
  \end{cases}.
\end{displaymath}

In this talk I take $\g=\affsl{2}$, so $I=\set{0,1}$ and
$c=h_0+h_1$.

\subsection{Perfect crystal base}
I recall briefly some facts from the theory of crystals bases
(see~\cite{K} for the definitions).  Recall that a \emph{crystal} is a
set~$B$ with maps $\te_i,\tf_i:B\sqcup\set{0} \rightarrow B\sqcup
\set{0}$ ($i\in I$), such that
\begin{enumerate}
\item $\te_i 0=0=\tf_i 0$,
\item there exists $n\in\Zplus$ such that $\te_i^n b=0=\tf_i^n b$
  ($\forall b\in B, i\in I$),
\item $b'=\tf_i b \iff b=\te_i b'$ ($b,b'\in B$, $i\in I$)
\end{enumerate}

Let $V$ be a finite dimensional $\U_q(\g')$ module with crystal base
$(L,B)$.  We take the $\U_q(\affsl{2}')$-module to be
$V:=\Directsum_{j\in J:=\range{0,2}}\Q(q)\v{j}$ with crystal
$B=\bigsqcup_{j\in J}\set{\b{j}}$ and crystal graph
\begin{displaymath}
  \b{0}\rightleftcrystal^1_0 \b{1}\rightleftcrystal^1_0 \b{2}.
\end{displaymath}
In the crystal graph an arrow $b\overset{i}{\rightarrow} b'
\iff b'=\tf_i b$.

Let $A:=\set{f\in\Q(q)\mid f\text{ has no pole at }q=0}$.  The
\emph{crystal lattice} $L= A\tensor_\Q \Span_\Q(B)$.

Define maps $\e_i,\f_i:B\rightarrow\N$ by
\begin{align*}
  \e_i(b)&:=\max\set{n\in\N\mid \te_i^n b\neq 0},
  \\
  \f_i(b)&:=\max\set{n\in\N\mid \tf_i^n b\neq 0}.
\end{align*}
and maps $\e,\f:B\rightarrow\Pcl$ by
\begin{displaymath}
  \e(b):=\sum_{i\in I}\varepsilon_i(b)\Lambda_i^{\cl},
  \qquad
  \f(b):=\sum_{i\in I}\varphi_i(b)\Lambda_i^{\cl}.
\end{displaymath}
In our example
\begin{equation}\label{epsilson-phi}
  \begin{aligned}
    \e(\b{j}) &= (2-j)\Lambda^{\cl}_0 + j\Lambda^{\cl}_1
    \\
    \f(\b{j}) &= j\Lambda^{\cl}_0 + (2-j)\Lambda^{\cl}_1
  \end{aligned}
  \qquad (j\in J).
\end{equation}
The weight of $b\in B$ is given by $\wt(b)=\f(b)-\e(b)$:
\begin{displaymath}
  \wt(\b{j})= 2(1-j)(\Lambda^{\cl}_1 -\Lambda^{\cl}_0).
\end{displaymath}

We take $\set{\v{j}=G(\b{j})}_{j\in J}$ to be a \emph{lower global
  base}~\cite{K} of $V$.  For our example, we have
\begin{align*}
  e_iG(b)&=[\f_i(b)+1]G(\te_i b),
  \\
  f_iG(b)&=[\e_i(b)+1]G(\tf_i b),
  \\
  q^hG(b)&=q^{\pairing{h,\wt(b)}} G(b).
\end{align*}
The use of a lower global base is essential to our construction: see
the remarks after Lemma~\ref{lem:condition} and
Theorem~\ref{thm:annihil}.

Let $B=\bigsqcup_{\lambda\in\Pcl}B_\lambda$ be the weight decomposition.

Let $B_1,B_2$ be crystals.  The tensor product $B_1\tensor B_2$
(corresponding to the coproduct $\Delta$) is defined to be the set
$B_1\times B_2$ with the action of the Kashiwara operators
$\te_i,\tf_i$ given by
\begin{align*}
  \te_i(b\tensor b') &=
  \begin{cases}
    \te_ib\tensor b' & \text{if $\varepsilon_i(b)>\varphi_i(b')$}\\
    b\tensor\te_ib' & \text{if $\varepsilon_i(b)\le\varphi_i(b')$}
  \end{cases}\: ,
  \\
  \tf_i(b\tensor b') &=
  \begin{cases}
    \tf_ib\tensor b'&\text{if $\varepsilon_i(b)\ge\varphi_i(b')$}\\
    b\tensor\tf_ib'&\text{if $\varepsilon_i(b)<\varphi_i(b')$}
  \end{cases}\: .
\end{align*}

The crystal base $(L,B)$ is perfect of level~$k=2$, which means that
it satisfies the following conditions:
\begin{enumerate}
\item There exists a weight $\lambda^\circ\in\wt(B)$ such that
  $\wt(B)$ is contained in the convex hull of $W\lambda^\circ$
  and that $\card{B_{w\lambda^\circ}}=1$ ($w\in W$).  The elements of
  $B_{w\lambda^\circ}$ and their weights $w\lambda^\circ$ are called
  \emph{extremal}.
\item $B\tensor B$ is connected as a crystal graph.
\item $k=\max\set{n\in\Zplus\mid \pairing{c,\e(b)}\geq n\: (\forall b\in
    B)}$.
\item $\e,\f:\set{b\in B\mid \pairing{c,\e(b)}=k}\rightarrow
  (\Pcl^+)_k:=\set{\lambda\in \Pcl^+\mid \pairing{c,\lambda}=k}$ are
  bijective.
\end{enumerate}
Note that $\pairing{c,\f(b)-\e(b)}=0$, since $\pairing{c,\wt(b)}=0$
($b\in B$).

In our example $\lambda^\circ=\pm 2(\Lambda^{\cl}_1 -\Lambda^{\cl}_0)$,
\begin{displaymath}
  B\tensor B=\quad
  \begin{CD}
    \b{0}\tensor\b{0} @>1>> \b{0}\tensor\b{1} @>1>> \b{0}\tensor\b{2}
    \\
    @AA0A  @AA0A @VV1V
    \\
    \b{1}\tensor\b{0} @>1>> \b{1}\tensor\b{1} @<0<< \b{1}\tensor\b{2}
    \\
    @AA0A @VV1V @VV1V
    \\
    \b{2}\tensor\b{0} @<0<< \b{2}\tensor\b{1} @<0<< \b{2}\tensor\b{2}
  \end{CD}
  \quad
\end{displaymath}
and \eqref{epsilson-phi} implies that $\e,\f$ map bijectively.
$\b{0},\b{2}$ are extremal.

\subsection{Affinization}
Let $\Vaff:=V\tensor\Q[z,z^{-1}]$ be the $\U_q(\g)$-module, which is
the affinization of $V$ such that
\begin{gather*}
  \begin{aligned}
    e_i (v\tensor\xi) &:= (e_i v)\tensor z^{\delta_{i,0}}\xi,
    \\
    f_i (v\tensor\xi) &:= (f_i v)\tensor z^{-\delta_{i,0}}\xi,
  \end{aligned}
  \quad (v\in V,\xi\in\Q[z,z^{-1}]),
  \\
  \wt(\v{j}\tensor z^a) := \wt(\v{j}) + a\delta, \qquad (j\in J, a\in\Z).
\end{gather*}
Let $(\Laff,\Baff)$ be the crystal base of $\Vaff$.

Let $z^a:\Vaff\rightarrow\Vaff$ ($a\in\Z$) denote the
$\U_q(\g')$-linear endomorphism $v\tensor\xi\mapsto v\tensor z^a\xi$
($v\tensor\xi\in\Vaff$).  Then $z^a\v{j}:=z^a(\v{j}\tensor
1)\equiv\v{j}\tensor z^a$ ($j\in J$, $a\in\Z$).

The crystal $\Baff=\bigsqcup_{j\in J,a\in\Z}\set{z^a\b{j}}$.  The
crystal graph of $\Baff$ has the following structure.
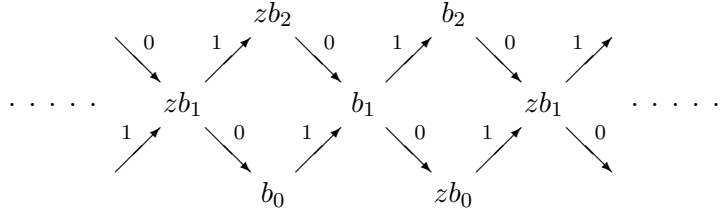
\begin{displaymath}
\setlength{\unitlength}{0.3mm}
\begin{picture}(300,89)(-10,735)
\put( 70,790){\vector( 1, 1){20}}\put( 70,770){\vector( 1,-1){20}}
\put(110,750){\vector( 1, 1){20}}\put(110,810){\vector( 1,-1){20}}
\put(150,790){\vector( 1, 1){20}}\put(150,770){\vector( 1,-1){20}}
\put(190,750){\vector( 1, 1){20}}\put(190,810){\vector( 1,-1){20}}
\put(230,790){\vector( 1, 1){20}}\put(230,770){\vector( 1,-1){20}}
\put( 30,810){\vector( 1,-1){20}}\put( 30,750){\vector( 1, 1){20}}
\multiput(260,780)(9,0){5}%
{\makebox(0.4444,0.6667){.}}
\multiput(20,780)(-9,0){5}%
{\makebox(0.4444,0.6667){.}}
\put( 60,775){\makebox(0,0)[b]{$z\b{1}$}}
\put(100,815){\makebox(0,0)[b]{$z\b{2}$}}
\put(100,735){\makebox(0,0)[b]{$\b{0}$}}
\put(140,775){\makebox(0,0)[b]{$\b{1}$}}
\put(180,815){\makebox(0,0)[b]{$\b{2}$}}
\put(180,735){\makebox(0,0)[b]{$z\b{0}$}}
\put(220,775){\makebox(0,0)[b]{$z\b{1}$}}
\put( 45,805){\makebox(0,0)[b]{\scriptsize 0}}
\put( 85,765){\makebox(0,0)[b]{\scriptsize 0}}
\put( 35,765){\makebox(0,0)[b]{\scriptsize 1}}
\put( 75,805){\makebox(0,0)[b]{\scriptsize 1}}
\put(115,765){\makebox(0,0)[b]{\scriptsize 1}}
\put(155,805){\makebox(0,0)[b]{\scriptsize 1}}
\put(195,765){\makebox(0,0)[b]{\scriptsize 1}}
\put(235,805){\makebox(0,0)[b]{\scriptsize 1}}
\put(125,805){\makebox(0,0)[b]{\scriptsize 0}}
\put(165,765){\makebox(0,0)[b]{\scriptsize 0}}
\put(205,805){\makebox(0,0)[b]{\scriptsize 0}}
\put(245,765){\makebox(0,0)[b]{\scriptsize 0}}
\end{picture}
\end{displaymath}
The arrows of the crystal graph induce a partial ordering of $\Baff$.

Define the map $l:\Baff\rightarrow\Z$ by
\begin{displaymath}
  l(z^a\b{j}):=2a-j.
\end{displaymath}

Define the energy function (\cite{KMN}) $H:\Baff\tensor\Baff\rightarrow\Z$ by
\begin{enumerate}
\item $H(zb\tensor b')=H(b\tensor b')-1$ and $H(b\tensor
  z b')=H(b\tensor b')+1$ ($b, b'\in\Baff$),
\item $H$ is constant on every connected component of the crystal
  graph $B_{\aff}\tensor B_{\aff}$,
\item $H(z^ab^\circ,z^ab^\circ)=0$ for any extremal $b^\circ\in B$ and
  $a\in\Z$.
\end{enumerate}
For our example we have
\begin{displaymath}
  H(\b{i}\tensor\b{j})=\min\set{i,2-j} \qquad (i,j\in J).
\end{displaymath}
If $H(b\tensor b')\leq 0$, then $l(b)\geq l(b')$.

\section{Wedge space}

For an extremal $b^\circ\in B$, let $v^\circ:=G(b^\circ)$.  Define
\begin{displaymath}
  N:= \U_q(\g)[z\tensor z,z^{-1}\tensor z^{-1},z\tensor 1+1\tensor z]
  (v^\circ\tensor v^\circ).
\end{displaymath}
It is not too difficult to prove that $N$ is independent of the choice
of extremal element $b^\circ$: indeed for any extremal $b'\in B$,
$G(z^a b')\tensor G(z^a b')\in N$ ($a\in\Z$).

For our example, define the following elements of $N$:
\begin{align*}
  C_{\b{0},\b{0}}&:= \v{0}\tensor\v{0},
  \\
  C_{\b{0},\b{1}}&:= \v{0}\tensor\v{1} +q^2\v{1}\tensor\v{0},
  \\
  C_{\b{0},\b{2}}&:= \v{0}\tensor\v{2} +q\v{1}\tensor\v{1}
  +q^4\v{2}\tensor\v{0},
  \\
  C_{\b{1},\b{2}}&:= \v{1}\tensor\v{2} +q^2\v{2}\tensor\v{1},
  \\
  C_{\b{2},\b{2}}&:= \v{2}\tensor\v{2},
  \\
  C_{z\b{2},\b{1}}&:= z\v{2}\tensor\v{1} +q^2\v{1}\tensor\v{2},
  \\
  C_{z^2\b{2},\b{0}}&:= z^2\v{2}\tensor\v{0} +qz\v{1}\tensor z\v{1}
  +q^4\v{0}\tensor z^2\v{2},
  \\
  C_{z\b{1},\b{0}}&:= z\v{1}\tensor\v{0} +q^2\v{0}\tensor z\v{1},
  \\
  C_{z\b{1},\b{1}}&:= z\v{1}\tensor\v{1} +q^2\v{1}\tensor z\v{1} +q^2
  [2](\v{0}\tensor z\v{2} + z\v{2}\tensor\v{0}).
\end{align*}
(These elements are constructed by the action of $\U_q(\affsl{2})$ on
$\v{1}\tensor\v{1}$.)
\begin{lemma}\label{lem:condition}
  Let $\mathcal{B}_{i,j}:=\set{(b,b')\in\Baff\tensor\Baff\mid
    H(b\tensor b')>0, l(\b{j})\leq
    l(b)<l(z^{H(\b{i}\tensor\b{j})}\b{i}), l(\b{j})<l(b')\leq
    l(z^{H(\b{i}\tensor\b{j})}\b{i})}$.  Each element
  $C_{z^{H(\b{i}\tensor\b{j})}\b{i},\b{j}}$ has the form
  \begin{displaymath}
    G(z^{H(\b{i}\tensor\b{j})}\b{i})\tensor G(\b{j})-
    \sum_{(b,b')\in\mathcal{B}_{i,j}} a_{b,b'}\,G(b)\tensor G(b')\, .
  \end{displaymath}
  The coefficients $a_{b,b'}$ lie in $q\,\Z[q]$.
  $H(z^{H(\b{i}\tensor\b{j})}\b{i}\tensor\b{j})= 0$.
\end{lemma}
The conditions in this Lemma are required for the construction of the
Fock space.
Note that this Lemma does not hold for the corresponding elements of
$N$ in an upper global base.

The vectors $\set{C_{z^{H(\b{i}\tensor\b{j})}\b{i},\b{j}}}$ are
linearly independent.
We have
\begin{displaymath}
  \sum_{i,j\in J}\Q(q)[z\tensor z,z^{-1}\tensor z^{-1},z\tensor
  1+1\tensor z] C_{z^{H(\b{i}\tensor\b{j})}\b{i},\b{j}} = N\,.
\end{displaymath}

Define the wedge product by
\begin{displaymath}
  \Wedge^2 \Vaff:= \Vaff\tensor\Vaff/N.
\end{displaymath}
For $v,v'\in\Vaff$, denote the image of $v\tensor
v'\in\Vaff^{\tensor2}$ in $\Wedge^2\Vaff$ by $v\wedge v'$.

We call a pair $(b,b')\in\Baff^{\tensor2}$ \emph{normally ordered},
if $H(b\tensor b')>0$.  For such a pair, $G(b)\wedge G(b')$ is
called a \emph{normally ordered wedge}.

Note that $\v{1}\wedge\v{1}$ is normally ordered and is \emph{not}
equal to 0 even at $q=1$.

The elements $C_{b,b'}\in N$ ($b,b'\in\Baff$
such that $H(b\tensor b')=0$) should be
thought of as a rule for writing $G(b)\wedge G(b')$ as a linear
combination of normally ordered wedges.

\begin{proposition}
  $\set{z^{a_1}\v{j_1}\wedge z^{a_2}\v{j_2}}_{H(z^{a_1}\b{j_1}\tensor
    z^{a_2}\b{j_2})>0}$ is a base of $\Wedge^2\Vaff$.
\end{proposition}

Let $n\in\Zplus$.  Define the $n$-wedge space by
\begin{align*}
  N_n &:=\sum_{r=0}^{n-2} \Vaff^{\tensor r}\tensor N\tensor
  \Vaff^{\tensor(n-r-2)},
  \\
  \Wedge^n\Vaff &:= \Vaff^{\tensor n}/N_n.
\end{align*}
By construction $\Wedge^n\Vaff$ is a $\U_q(\g)$-module.

A sequence $(\b{j_1},\b{j_2},\dots,\b{j_n})\in \Baff^{\tensor n}$
($j_r\in J$) is called \emph{normally ordered} if
$H(\b{j_m}\tensor\b{j_{m+1}})>0$ for every $m\in\range{0,n-1}$.  In
this case the corresponding wedge $\v{j_1}\wedge\v{j_2}\wedge\cdots
\wedge\v{j_n}$ is called a \emph{normally ordered wedge}.

The homomorphism
\begin{align*}
  \Wedge^{r_1}\Vaff\tensor\Wedge^{r_2}\Vaff &\rightarrow
  \Wedge^{r_1+r_2}\Vaff
  \\
  u_1\tensor u_2 &\mapsto u_1\wedge u_2
\end{align*}
is $\U_q(\g)$-linear.

\begin{theorem}
  Normally ordered $n$-wedges form a base of $\Wedge^n \Vaff$.
\end{theorem}

Define $L(\Wedge^n\Vaff)$ to be the image of $L(\Vaff^{\tensor n})$ in
$\Wedge^n \Vaff$ and
$B(\Wedge^n\Vaff):=\set{(\b{j_1},\b{j_2},\dots,\b{j_n})\mid
  \text{normally ordered}}$.

\section{Level 2 Fock space}

\subsection{Ground state sequence}
Recall that $\Baff$ is the affinization of the perfect level~$k=2$
crystal $B$.  Extend the maps $\e,\f:B\rightarrow \Pcl$ to
$\Baff\rightarrow\Pcl$ by $\e(z^a\b{j}):=\e(\b{j})$ and
$\f(z^a\b{j})=\f(\b{j})$.  Fix a sequence $\bo{m}\in\Baff$ ($m\in\Z$)
such that
\begin{gather*}
  \begin{aligned}
    \pairing{c,\e(\bo{m})}&=k,
    \\
    \e(\bo{m})&=\f(\bo{m+1}),
  \end{aligned}
  \\
  H(\bo{m}\tensor\bo{m+1})=1.
\end{gather*}
We call this sequence a \emph{ground state sequence}.  Take also a
sequence of weights $\lambda_m$ ($m\in\Z$) such that
\begin{gather*}
  \begin{aligned}
    \pairing{c,\lambda_m} &= k,
    \\
    \cl(\lambda_m)&=\varphi(b^\circ_{m}),
  \end{aligned}
  \\
  \lambda_m = \wt(\bo{m}) +\lambda_{m+1}.
\end{gather*}
Define $\vo{m}:=G(\bo{m})$.

For our example of level~2 $\U_q(\affsl{2})$ there are, up to
equivalence, only two possible ground state sequences
\begin{equation}\tag{A}
  \begin{aligned}
    \bo{m} &=
    \begin{cases}
      z\b{2} &\text{for } m\in 2\Z,
      \\
      \b{0} &\text{for } m\in2\Z+1,
    \end{cases}
    \\
    \lambda_m &=
    \begin{cases}
      2\Lambda_0 +\delta &\text{for } m\in 2\Z,
      \\
      2\Lambda_1 &\text{for } m\in2\Z+1,
    \end{cases}
  \end{aligned}
\end{equation}
and
\begin{equation}\tag{B}
  \begin{aligned}
    \bo{m} &= \b{1} \qquad\qquad (m\in\Z),
    \\
    \lambda_m &= \Lambda_0 +\Lambda_1 \qquad (m\in\Z).
  \end{aligned}
\end{equation}

\subsection{Fock space}

In the formal semi-infinite wedge space $\Wedge^\infty\Vaff$, define
the vacuum vector
\begin{displaymath}
  \overline{\vac{m}}:= \vo{m}\wedge\vo{m+1}\wedge\vo{m+2}\wedge\cdots
  \qquad (m\in\Z).
\end{displaymath}
Define the pre-Fock space
\begin{displaymath}
  \bar\F{m}:= \sum_{r\in\N} \Wedge^r\Vaff\wedge\overline{\vac{m+r}}.
\end{displaymath}
Let $L(\bar\F{m})$ be the crystal lattice of $\bar\F{m}$.
Define the Fock space $\F{m}$ to be
\begin{displaymath}
  \F{m}:= \bar\F{m}/ \bigcap_{n>0} q^n L(\bar\F{m}).
\end{displaymath}
Let $\vac{m}$ be the image of $\overline{\vac{m}}$ in $\F{m}$.  Taking
$\set{q^nL(\F{m})}_{n\in\N}$ as a neighborhood system of $0$, $\F{m}$
is endowed with a $q$-adic topology.  By construction the $q$-adic
topology is separated, since $\bigcap_{n>0} q^n L(\F{m})=0$.

We have an algebra homomorphism $\Wedge^r\Vaff\tensor
\F{m}\rightarrow\F{m-r}$.

\newcommand{\tA}{{\text{A}}}
\newcommand{\tB}{{\text{B}}}

Denote the vacuum vector $\vac{m}$ associated to ground sequences (A)
and (B) by $\vac{m}^\tA$ and $\vac{m}^\tB$ respectively.  Similarly
denote the Fock space $\F{m}$ associated to ground sequences (A) and
(B) by $\F{m}^\tA$ and $\F{m}^\tB$ respectively.

Note that $G(b)\wedge\overline{\vac{m}}=0$, if and only if there exists
$r\in\N$ such that
$G(b)\wedge\vo{m}\wedge\vo{m+1}\wedge\cdots\wedge\vo{m+r}=0$.
\begin{theorem}\label{thm:annihil}
  Let $b\in \Baff$ be such that $H(b\tensor\bo{m})\leq 0$. In $\F{m}$ the
  equality
  \begin{displaymath}
    G(b)\wedge\vac{m}=0
  \end{displaymath}
  holds in the $q$-adic topology.
\end{theorem}
Note that this theorem only holds for a lower global base.

For example in $\F{m}^\tB$ the vacuum vector is
$\vac{m}^\tB=\v{1}\wedge\v{1}\wedge\v{1}\wedge\dots$ and
\begin{equation}\label{vacuum-cond}
  \begin{aligned}
    \v{0}\wedge\vac{m}^\tB &= \lim_{r\rightarrow\infty}(-q^2)^r
    \v{1}^{\wedge r}\wedge\v{0}\wedge\vac{m+r}^\tB
    \\
    &= 0 \;\text{ (in the $q$-adic topology)}.
  \end{aligned}
\end{equation}

For a normally ordered sequence
$(z^{a_m}\b{j_m},z^{a_{m+1}}\b{j_{m+1}},\dots)$ in $\Baff$ such that
$z^{a_r}\b{j_r}=\bo{r}$ for all $r\gg m$, the wedge
$z^{a_m}\v{j_m}\wedge z^{a_{m+1}}\v{j_{m+1}}\wedge\cdots$ is called a
\emph{normally ordered wedge}.

\begin{theorem}
  The normally ordered wedges in $\F{m}$ form a base of $\F{m}$.
\end{theorem}

\subsection{$\U_q(\g)$-module structure}

First we assign weights to the Fock space by setting
\begin{displaymath}
  \wt(\vac{m}):= \lambda_m.
\end{displaymath}
This fixes the action of the $q^h$ ($h\in P^*$).
For example $\wt(\vac{m}^\tB)= \Lambda_0 +\Lambda_1$ $(m\in\Z)$.

\begin{proposition}
  Let $V(\lambda_m)$ denote the irreducible integrable $\U_q(\g)$-module of
  highest weight $\lambda_m$.  The character of $\F{m}$ is
  \begin{equation}\label{character}
    \ch(\F{m})=\ch(V(\lambda_m))\prod_{r>0}(1-\mathrm{e}^{-r\delta})^{-1}.
  \end{equation}
\end{proposition}

\begin{theorem}
  $\F{m}$ has the structure of an integrable $\U_q(\g)$-module.
\end{theorem}

There are two steps to proving this: (i)~proving that the action of
$e_i,f_i$ is well-defined (converges in the $q$-adic topology)
and is integrable, and (ii) checking that the commutation relations are
satisfied.  In fact it is sufficient to prove it on the vacuum vector
$\vac{m}$.

For example in $\F{m}^\tB$ we have
\begin{align*}
  t_i \,\vac{m}^\tB &= q\,\vac{m}^\tB \qquad (i\in I),
  \\
  e_1 \,\vac{m}^\tB &= [2]\sum_{r\in\N} \v{1}^{\wedge r}\wedge
  \v{0}\wedge
  \vac{m+1+r}^\tB
  = 0\; \text{ (by~\eqref{vacuum-cond})},
  \\
  f_1 \,\vac{m} &= [2]\sum_{r\in\N} \v{1}^{\wedge r}\wedge \v{2}\wedge
  q\,
  \vac{m+1+r}^\tB\\
  &= q\,[2]\Bigl(\sum_{r\in\N} (-q^2)^r\Bigr)\, \v{2}\wedge
  \vac{m+1}^\tB = \v{2}\wedge \vac{m+1}^\tB.
\end{align*}
Then one can check that
\begin{align*}
  [e_1,f_1] \,\vac{m}^\tB = e_1\cdot f_1\,\vac{m}^\tB
  &= e_1 \v{2}\wedge \vac{m+1}^\tB
  \\
  &= \vac{m}^\tB = \frac{t_1-t_1^{-1}}{q-q^{-1}}\,\vac{m}^\tB.
\end{align*}

\subsection{Bosons}
Define boson operators $B_a$
\begin{align*}
  B_a :=& \sum_{r\in\N} 1^{\tensor r}\tensor z^a \tensor
  1^{\tensor\infty} \qquad (a\in\Z\setminus\set{0})
  \\
  &\equiv z^a \tensor 1^{\tensor\infty} +1\tensor z^a \tensor 1^{\tensor\infty}
  +1\tensor1\tensor z^a \tensor 1^{\tensor\infty} +\dots \; .
\end{align*}

\begin{proposition}
  The action of the operators $B_a$ on $\F{m}$ converges in the
  $q$-adic topology.
\end{proposition}

\begin{proposition}
  $B_a\,\vac{m}=0$ for all $a\in\Zplus$.
\end{proposition}

\begin{proposition}
  There exists $\gamma_a\in\Q(q)$ (independent of $m$) such that
  \begin{displaymath}
    [B_a,B_{a'}]\:\vac{m}= \gamma_a \,\delta_{a+a',0}\,\vac{m}.
  \end{displaymath}
  At $q=0$, $\gamma_a=a$.
\end{proposition}

Let $\H$ be the Heisenberg algebra generated by
$\set{B_a}_{a\in\Z\setminus\set{0}}$ with the defining relations
$[B_a,B_{a'}]=\delta_{a+a',0}\gamma_a$.  Then $\H$ acts on the Fock
space $\F{m}$ commuting with the action of $\U_q(\g')$.  Let
$\Q[\H_-]:=\Q[B_{-a}]_{a\in\Zplus}\cdot 1$ be the Fock space for $\H$
with vacuum vector~$1$ and the defining relation $B_a\cdot1=0$
($a\in\Zplus$).  Let $u_{\lambda_m}$ denote the highest weight vector
in $V(\lambda_m)$.  Since $\vac{m}$ is annihilated by $e_i$ ($i\in I$)
and $B_a$ ($a\in\Zplus$), we have an injective $\U_q(\g')\tensor
\H$-linear homomorphism
\begin{align*}
  \iota_m:V(\lambda_m)\tensor\Q[\H_-] &\rightarrow \F{m}
  \\
  u_{\lambda_m}\tensor 1 &\mapsto \vac{m}.
\end{align*}

\begin{theorem}
  $\iota_m:V(\lambda_m)\tensor\Q[\H_-]\rightarrow \F{m}$ is an
  isomorphism.
\end{theorem}
The proof is by comparing the characters~\eqref{character}.  This
gives the decomposition of $\F{m}$ into irreducible
$\U_q(\g)$-modules.

$\gamma_a$ can be calculated by using the decomposition via $\iota_m$
of the wedge vertex operator ($\Vaff\tensor\F{m}\rightarrow\F{m-1}$)
to a product of the usual vertex operator ($\Vaff\tensor
V(\lambda_m)\rightarrow V(\lambda_{m-1})$) and the boson vertex
operator, and then calculating the equality of two-point functions
corresponding to this decomposition (see~\cite{KMPY}).
\begin{proposition}
  In our example
  \begin{displaymath}
    [B_a,B_{a'}]\,\vac{m}^\tB= \frac{a}{1-q^{2a}}
    \,\delta_{a+a',0}\,\vac{m}^\tB.
  \end{displaymath}
\end{proposition}

\medskip
\begin{flushleft}\small
  R.I.M.S., Kyoto University, Kyoto 606-01, Japan
  \\
  \verb|petersen@kurims.kyoto-u.ac.jp|
  \\
  \verb|http://www.kurims.kyoto-u.ac.jp/~petersen/|
\end{flushleft}
\end{document}